\documentstyle{amsppt}
\NoRunningHeads
\magnification=\magstep1
\hyphenation{co-deter-min-ant co-deter-min-ants pa-ra-met-rised
pre-print
fel-low-ship Cox-et-er dis-trib-ut-ive}
\def\leaderfill{\leaders\hbox to 1em{\hss.\hss}\hfill}
\def\A{{\Cal A}}

\def\H{{\Cal H}}

\def\idest{i.e.\ }
\def\a{{\alpha}}
\def\b{{\beta}}
\def\g{{\gamma}}
\def\G{{\Gamma}}
\def\d{{\delta}}
\def\e{{\varepsilon}}
\def\z{{\zeta}}

\def\l{{\lambda}}

\def\dt{{\Bbb D \Bbb T}}
\def\boxit#1{\vbox{\hrule\hbox{\vrule \kern3pt
\vbox{\kern3pt\hbox{#1}\kern3pt}\kern3pt\vrule}\hrule}}
\def\rabbit{\vbox{\hbox{\kern0pt
\vbox{\kern0pt{\hbox{---}}\kern3.5pt}}}}

\def\tableau#1{
        \hbox {
                \hskip -10pt plus0pt minus0pt
                \raise\baselineskip\hbox{
                \offinterlineskip
                \hbox{#1}}
                \hskip0.25em
        }
}

\def\tabCol#1{
\hbox{\vtop{\hrule
\halign{\strut\vrule\hskip0.5em##\hskip0.5em\hfill\vrule\cr\lower0pt
\hbox\bgroup$#1$\egroup \cr}
\hrule
} } \hskip -10.5pt plus0pt minus0pt}

\def\CR{
        $\egroup\cr
        \noalign{\hrule}
        \lower0pt\hbox\bgroup$
}



\def\xfanb{{\bf 1}}

\def\xfga{{\bf 2}}
\def\xfy{{\bf 3}}
\def\xgra{{\bf 4}}
\def\xgl{{\bf 5}}

\def\xrmgl{{\bf 6}}

\def\xhum{{\bf 7}}

\def\xjon{{\bf 8}}

\def\xk{{\bf 9}}

\def\xmsa{{\bf 10}}

\def\xsteb{{\bf 11}}
\def\xtl{{\bf 12}}
\def\xwes{{\bf 13}}

\topmatter
\title
Cellular algebras arising from Hecke algebras of type $H_n$
\endtitle

\author R.M. Green\endauthor
\affil 
Mathematical Institute\\ Oxford University\\ 24--29 St. Giles'\\
Oxford OX1 3LB\\ England\\ 
{\it  E-mail:} greenr\@maths.ox.ac.uk
\endaffil

\abstract
We study a finite-dimensional quotient of the Hecke algebra of type
$H_n$ for general $n$, using a calculus of diagrams.  This provides a
basis of monomials in a certain set of generators.  Using this, we
prove a conjecture of C.K. Fan about the semisimplicity of the
quotient algebra.  We also discuss the cellular structure of the
algebra, with certain restrictions on the ground ring.
\endabstract


\thanks
The author was supported in part by an E.P.S.R.C. postdoctoral
research assistantship.
\endthanks
\endtopmatter

\centerline{\bf To appear in ``Mathematische Zeitschrift''}

\head 0. Introduction \endhead

There has been much recent interest in the Temperley--Lieb algebra and
its various generalisations.  Graham [\xgra] in his thesis studied a
certain quotient, which we will call $TL(X)$,
of a Hecke algebra $\H(X)$ associated to a Dynkin
diagram $X$.  In the case where $X$ is a Dynkin diagram of type $A$,
this quotient was considered by Jones [\xjon], who pointed out that it is
nothing other than the Temperley--Lieb algebra, which first appeared
in [\xtl].  The Temperley--Lieb algebra has applications in several areas of
mathematics, including statistical mechanics and knot theory.

A remarkable feature of the algebras $TL(X)$ is that they can be
finite dimensional, even when $\H(X)$ is infinite dimensional.  Graham
[\xgra] classified the finite dimensional algebras $TL(X)$
into seven infinite families: $A, B, D, E, F, H$ and $I$.  (Contrast
this to the classification of Hecke algebras associated to irreducible
Coxeter systems, in which there are only
finitely many algebras of types $E$, $F$ and $H$.)  

This paper is concerned with the infinite family of type $H$, in which case
the Hecke algebra
$\H(H_n)$ is finite dimensional only for $n \leq 4$.  The algebra $TL(H_n)$ was
mentioned briefly by Fan in [\xfanb, \S7.3], where it was conjectured that
$TL(H_n)$ is generically semisimple.  The dimensions of the
generically irreducible modules are also conjectured.
In the course of the paper, we will prove these conjectures.
Note that semisimplicity is obvious in the cases
where $\H(X)$ is finite dimensional, because in this case $\H(X)$ is
itself generically semisimple, but this argument fails if $\H(X)$ is
infinite dimensional. 

Our approach is first to realise $TL(H_n)$ as an algebra of diagrams
arising from the category of ``decorated tangles'' which was
introduced by the author in [\xrmgl].  Diagram calculi have
already been developed for algebras of some of the other types: for $TL(A_n)$
it is well known (see [\xwes] or [\xgl, \S6]), types $B_n$ and $D_n$
were done in [\xrmgl], and the infinite-dimensional
``affine'' Temperley--Lieb algebra $TL(\widehat A_n)$ was tackled in
[\xfga].  This is interesting to do in its own right, since many
natural questions about
the algebras (such as the determination of the cells, dimensions and
structure constants) have simple formulations in terms of the
combinatorics of the associated diagrams.  We
will show that the diagrams may be adapted into a datum for a cellular
algebra (in the sense of [\xgl]), provided that the polynomial $x^2 -
x - 1$ splits into distinct linear factors over the ground ring.
Since the algebra $TL(H_2)$ is a $q$-analogue of a 9-dimensional
quotient of the group algebra of the dihedral group of order 10, it
seems that such a hypothesis cannot be usefully weakened.

The algebras $TL(X)$ of types $A, B, D, E$ and $F$ each have a basis
consisting of monomials in the obvious set of algebra generators
(which correspond to the Coxeter generators), and the
structure constants with respect to this basis are positive in a
natural sense.  Furthermore, the product of two monomials is a scalar
multiple of another monomial.
In type $H$, the obvious basis of monomials does not have the
positivity property, and it is not true that the product of two
monomials is a scalar multiple of one other.  This means that Fan's
techniques from [\xfanb] are unsuitable for analysing the algebra $TL(H_n)$.

In this paper, we overcome this problem by working with the basis of
diagrams, which has much more convenient properties (e.g. positivity of
structure constants and compatibility with cellular algebras).  This
basis is not obvious from the description of $TL(H_n)$ via generators and
relations, but is very natural from the viewpoint of decorated
tangles.  We also show how the new basis elements can be expressed as
monomials in a slightly larger set of algebra generators.

\vskip 20pt
\head 1. Preliminaries \endhead

\subhead 1.1 Coxeter groups of type $H_n$ \endsubhead

Let $n \in {\Bbb N}$ be at least $2$.
The Coxeter group of type $H_n$ corresponds to the Coxeter graph shown in
Figure 1.

\topcaption{Figure 1} Coxeter graph of type $H_n$ \endcaption
\centerline{
\hbox to 3.138in{
\vbox to 0.402in{\vfill
        \includegraphics{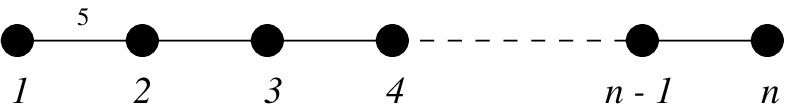}
}
\hfill}
}

\definition{Definition 1.1.1}
The Coxeter group $W(H_n)$ is given by generating involutions $\{s_i :
i \leq n\}$ and defining relations $$\eqalign{
s_i s_j &= s_j s_i \text{\quad if $|i - j| > 1$},\cr
s_i s_j s_i &= s_j s_i s_j \text{\quad if $|i - j| = 1$ and $\{i, j\}
\ne \{1, 2\}$},\cr
s_1 s_2 s_1 s_2 s_1 &= s_2 s_1 s_2 s_1 s_2.\cr
}$$
\enddefinition

\remark{Remark 1.1.2}
The group $W(H_2)$ is isomorphic to the dihedral group of order 10
(\idest $W(I_2)$), but it will be convenient to regard it as a Coxeter
group of type $H$ in some of our proofs.

As explained in [\xhum, \S2], the groups $W(H_n)$ are finite for $n =
2, 3, 4$, where they have orders $10$, $120$ and $14400$
respectively.  These groups occur as the full symmetry groups of Platonic
solids with pentagonal faces: $H_2$ corresponds to the pentagon, $H_3$
to the dodecagon and $H_4$ to a regular 120-sided solid in 4
dimensions.  For $n > 4$, the group $W(H_n)$ is infinite, which is
reminiscent of the fact that there is no analogue of the dodecagon in
higher dimensions, the only Platonic solids being generalized
tetrahedra, cubes and octahedra.
\endremark

\vskip 20pt

\subhead 1.2 Hecke algebras of type $H_n$ \endsubhead

We now introduce the Hecke algebra and its quotient $TL(H_n)$.

\definition{Definition 1.2.1}
The Hecke algebra $\H(H_n)$ is defined over the ring $$\A := {\Bbb
Z}[v, v^{-1}],$$ where $v = q^{1/2}$.  It has a free $\A$-basis $\{T_w
: w \in W(H_n)\}$, and the multiplication is defined by the rules $$
T_s T_w := \cases
T_{sw} & \text{ if } \ell(sw) > \ell(w),\cr
q T_{sw} + (q-1) T_w & \text{ otherwise.}\cr
\endcases$$  Here, $\ell(w)$ is the length of $w$, \idest the length
of a shortest word in the $s_i$ which is equal to $w$.
\enddefinition

Following Fan [\xfanb, \S7.3] and Graham [\xgra], 
we make the following definition.

\definition{Definition 1.2.2}
Let $n \in {\Bbb N} \geq 2$.  We define the associative, unital algebra
$TL(H_n)$ over $\A$ via generators $E_1, E_2, \ldots E_n$ and 
relations $$\eqalign{
E_i^2 &= [2] E_i, \cr
E_i E_j &= E_j E_i \text{\quad if \ $|i - j| > 1$},\cr
E_i E_j E_i &= E_i \text{\quad if \ $|i - j| = 1$ \ and \ $i, j > 1$},\cr
E_i E_j E_i E_j E_i &= 3 E_i E_j E_i - E_i
\text{\quad if \ $\{i, j\} = \{1, 2\}$}.\cr
}$$  Here, $[2]$ denotes the Laurent polynomial $v + v^{-1}$.
\enddefinition

\remark{Remark 1.2.3}
The algebra $TL(H_n)$ is a quotient of $\H(H_n)$
which corresponds to the Coxeter graph in Figure 1.  The quotient map
takes the Kazhdan--Lusztig basis element $C'_s = v^{-1}T_e +
v^{-1}T_s$, where $e$ is the identity and $s$ is of length 1, to
$E_s$.

Later, we shall want to replace the base ring $\A$ with a field, but
we will not be concerned with trying to generalise the results to 
characteristics 2, 3 or 5.
\endremark

It is convenient for later purposes to define the following elements
of $TL(H_n)$.

\definition{Definition 1.2.4}
We define $$\eqalign{
\a &:= E_1 E_2 - 1,\cr
\b &:= E_2 E_1 - 1,\cr
\e &:= E_1 E_2 E_1 - 2 E_1 \text{ and }\cr
\z &:= E_2 E_1 E_2 - 2 E_2.\cr
}$$\enddefinition

\remark{Remark 1.2.5}
Notice that we can rephrase the non-monomial relations in Definition
1.2.2 as the monomial relations $\e\b = E_1$ and $\z\a = E_2$.
\endremark

\vskip 20pt

\head 2. The diagram algebra $\Delta_n$ \endhead

We define a calculus of diagrams which will be seen in \S3 to describe
the generalized Temperley--Lieb algebra of type $H$.  A convenient way
to explain this is via the category of ``decorated tangles'' which was
introduced by the author in [\xrmgl].

\subhead 2.1 The category of decorated tangles \endsubhead

Following [\xfy], we define a tangle as follows.

\definition{Definition 2.1.1}
A tangle is a portion of a knot diagram contained in a rectangle.  The
tangle is incident with the boundary of the rectangle only on the
north and south faces, where it intersects transversely.  The
intersections in the north (respectively, south) face are numbered
consecutively starting with node number $1$ at the western (\idest the
leftmost) end.

Two tangles are equal if there exists an isotopy of the plane carrying
one to the other such that the corresponding faces of the rectangle
are preserved setwise.
\enddefinition

We call the edges of the rectangular frame ``faces'' to avoid
confusion with the ``edges'' which are the arcs of the tangle.

We extend the notion of a tangle so that each arc of the tangle may be
assigned a nonnegative
integer.  (This is similar to the notion of ``coloured'' tangles in
[\xfy].)  If an arc is assigned the value $r$, we represent this
pictorially by decorating the arc with $r$ blobs.  We also require
some further restrictions, as explained in the following definition.

\definition{Definition 2.1.2}
A decorated tangle is a crossing-free tangle in which each arc is
assigned a nonnegative integer.  Any arc not exposed to the west face
of the rectangular frame must be assigned the integer $0$.
\enddefinition

\remark{Remark 2.1.3}
This means that any decorated tangle consists only of loops and edges,
none of which intersect each other.
\endremark

\example{Example 2.1.4}
Figure 2 shows a typical example of a decorated tangle.  We will
tend to emphasise the intersections of the tangle with the frame
rather than the frame itself, which is why each node (\idest
intersection point with the frame) is denoted by a disc.  In this
case, the only edges or loops
exposed to the west wall are the three which already carry decorations.
\endexample

\topcaption{Figure 2} A decorated tangle \endcaption
\centerline{
\hbox to 3.027in{
\vbox to 0.888in{\vfill
        \includegraphics{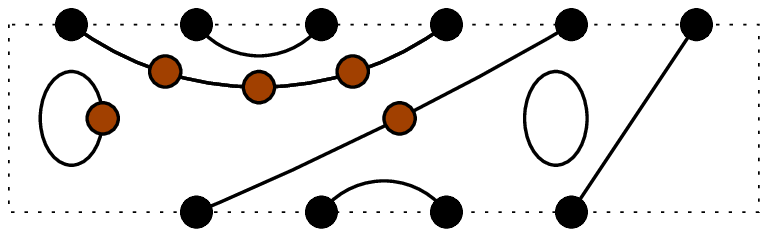}
}
\hfill}
}

We now define a category based on the set of decorated tangles, as
follows.

\definition{Definition 2.1.5}
The category of decorated tangles, $\dt$, has as its objects the
natural numbers (not including zero).  The morphisms from $n$ to $m$
are the decorated tangles with $n$ nodes in the north face and $m$ in
the south.  The
source of a morphism is the number of points in the north face of the
bounding rectangle, and the target is the number of points in the
south face.  Composition of morphisms works by concatenation of the
tangles, matching the relevant south and north faces together.
\enddefinition

\remark{Remark 2.1.6}
Note that for there to be any morphisms from $n$ to $m$, it is
necessary that $n+m$ be even.  Also notice that the asymmetric
properties of the west face of the rectangle mean that we cannot
introduce the tensor product of two morphisms by the lateral
juxtaposition of diagrams as in [\xfy].
\endremark

The category-theoretic definition allows us to define an algebra of
decorated tangles, as follows.

\definition{Definition 2.1.7}
Let $R$ be a commutative ring and let $n$ be a positive integer.  
Then the $R$-algebra $\dt_n$ has as a
free $R$-basis the morphisms from $n$ to $n$, where the multiplication
is given by the composition in $\dt$.
\enddefinition

\definition{Definition 2.1.8}
The edges in a tangle $T$ which connect nodes (\idest not the loops)
may be classified
into two kinds: propagating edges, which link a node in the north
face with a node in the south face, and non-propagating edges, which
link two nodes in the north face or two nodes in the south face.
\enddefinition

\vskip 20pt

\subhead 2.2 $H$-admissible diagrams \endsubhead

We introduce the concept of an $H$-admissible diagram, which plays a
key r\^ole in describing the diagram calculus relevant for $TL(H_{n-1})$.

\definition{Definition 2.2.1}
An $H$-admissible diagram with $n$ strands is an element of $\dt_n$
with no loops which satisfies the following conditions.

\item{\rm (i)}
{No edge may be decorated if all the edges are propagating.}
\item{\rm (ii)}
{If there are non-propagating edges in the diagram, then either there
is a decorated edge in the north face connecting nodes 1 and 2, or
there is a non-decorated edge in the north face connecting nodes $i$
and $i+1$ for $i > 1$.
A similar condition holds for the south face.}
\item{\rm(iii)}
{Each edge carries at most one decoration.}
\enddefinition

An example of an $H$-admissible diagram for $n = 6$ is shown in Figure
3.

\topcaption{Figure 3} An $H$-admissible diagram \endcaption
\centerline{
\hbox to 2.638in{
\vbox to 0.888in{\vfill
        \includegraphics{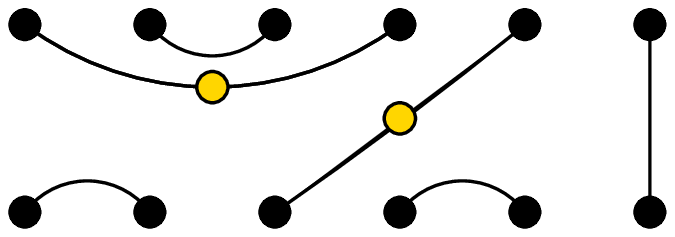}
}
\hfill}
}

The point of part (ii) in Definition 2.2.1 excludes situations like
the one where Figure 4 appears as the top half of an element of
$\dt_n$.

\topcaption{Figure 4} The north face of a diagram excluded by
Definition 2.2.1 (ii) \endcaption
\centerline{
\hbox to 3.138in{
\vbox to 0.277in{\vfill
        \includegraphics{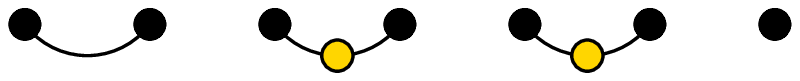}
}
\hfill}
}

\definition{Definition 2.2.2}
The algebra $\Delta_n$ (over a commutative ring with identity)
has as a basis the $H$-admissible
diagrams with $n$ strands and multiplication induced from that of
$\dt_n$ subject to the relations shown in Figure 5.
\enddefinition

\topinsert
\topcaption{Figure 5} Reduction rules for $\Delta_n$ \endcaption
\centerline{
\hbox to 1.041in{
\vbox to 1.916in{\vfill
        \includegraphics{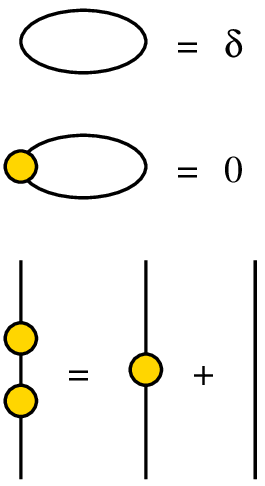}
}
\hfill}
}
\endinsert

\remark{Remark 2.2.3}
The meaning of the first relation in Figure 5 is that any
undecorated loop can be removed and the resulting tangle multiplied by
$\d$.  The second relation means that any tangle containing a loop
with one decoration is equivalent to $0 \in \Delta_n$.  The third
relation means that any tangle $T$ containing an edge or loop $\e$ with $r$ 
decorations $r > 1$ is equivalent to the sum of two other tangles $T'$
and $T''$ which are the same as $T$ except that the edge
or loop corresponding to $\e$ carries $r-1$ (respectively, $r-2$) decorations.

The second rule can be modified so that its removal corresponds to
multiplication by a second parameter, $\d'$.  This would eventually lead to a
two-parameter version of $TL(H_{n-1})$, but we do not pursue this here.
\endremark

\proclaim{Lemma 2.2.4}
The relations in Figure 5 allow the product of two elements of
$\Delta_n$ to be expressed unambiguously as a linear combination of
basis elements.  This makes $\Delta_n$ into an associative algebra.
\endproclaim

\demo{Proof}
We observe that the product of two $H$-admissible diagrams can be
expressed as a linear combination of others by using the reduction
rules given.

A case by case check shows that the order in which the relations are
applied is immaterial and that the end result can therefore be
expressed unambiguously in terms of the basis of $H$-admissible
diagrams.

Using these observations, associativity is inherited from the
associativity of $\dt_n$, by consideration of the concatenation of
three tangles $T T' T''$.
\qed\enddemo

\vskip 20pt

\head 3. Realisation of $TL(H_n)$ as an algebra of diagrams \endhead

\subhead 3.1 Representation of $TL(H_n)$ by diagrams \endsubhead

One of our main aims will be to show that $\Delta_{n+1}$ and $TL(H_n)$
are isomorphic.  To do this, we show how to represent $TL(H_n)$ using
the $H$-admissible diagrams.

\definition{Definition 3.1.1}
The $H$-admissible diagram $U_i$, where $1 \leq i \leq n$, is the diagram
all of whose edges are propagating and undecorated, except for those
attached to nodes $i$ and $i+1$ in the north row, and nodes $i$ and
$i+1$ in the south row.  These four nodes are connected in the pairs
given, using decorated edges if $i = 1$, and using undecorated edges
if $i > 1$.
\enddefinition

\demo{Examples}
When $n = 6$, the elements $U_1$ and $U_2$ are as shown in Figures 6
and 7.
\enddemo

\topcaption{Figure 6} The diagram $U_2$ for $n = 6$ \endcaption
\centerline{
\hbox to 2.638in{
\vbox to 0.888in{\vfill
        \includegraphics{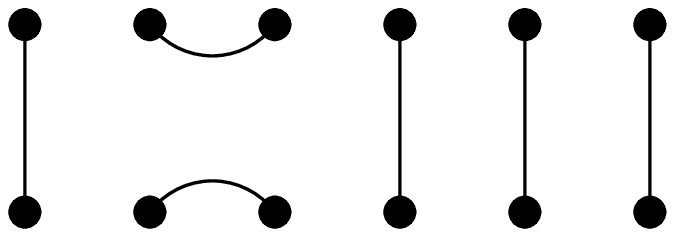}
}
\hfill}
}

\topcaption{Figure 7} The diagram $U_1$ for $n = 6$ \endcaption
\centerline{
\hbox to 2.638in{
\vbox to 0.888in{\vfill
        \includegraphics{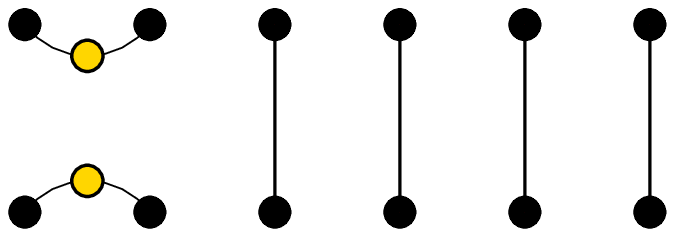}
}
\hfill}
}

From now on, we take the base ring for $\Delta_{n+1}$ to be $\A$, meaning
that the parameter $\d$ is $[2]$.  More general results may be found
by tensoring over a suitable ring.

\proclaim{Proposition 3.1.2}
There is a homomorphism of $\A$-algebras from $TL(H_n)$ to
$\Delta_{n + 1}$ which takes $E_i$ to $U_i$ for each $i$.
\endproclaim

\demo{Proof}
This is simply a matter of checking that all the relations in
Definition 1.2.2 hold, which presents no problems.
\qed\enddemo

\vskip 20pt

\subhead 3.2 Algebra generators for $\Delta_{n+1}$ \endsubhead

In order to prove that $\rho$ is an isomorphism, we
will first show that $\Delta_{n+1}$ is generated as an $\A$-algebra
(with identity) by the elements $U_i$.

During the course of the proofs, it helps to understand the
case $n = 2$, which was the motivation for Definition 1.2.4.

\proclaim{Lemma 3.2.1}
The map $\rho$ is an isomorphism for $n = 2$.  The basis of $H$-admissible
diagrams consists of the images of the $9$ elements $$
1, E_1, E_2, E_1 E_2, E_2 E_1, E_1 \b, E_2 \a, E_1 \z, E_2 \e = \z E_1.
$$  Thus, $\Delta_3$ is generated as an $\A$-algebra with $1$
by $U_1$ and $U_2$.
\endproclaim

\demo{Proof}
This is another routine exercise using the diagram multiplication,
which is instructive to carry out.
\qed\enddemo

To deal with the case for general $n$, it is convenient to introduce a
number of ``moves'', in which a diagram element is multiplied (on the
left or on the right) by a
monomial in the generators $$
G_n := \{1, U_1, \ldots, U_n, \rho(\a), \rho(\b), \rho(\z)\}
$$ to form another diagram element.  It will eventually turn out that
any $H$-admissible diagram may
be obtained as a suitable word in the generators $G_n$.

In the next five lemmas, $D$ is an $H$-admissible diagram.  The proofs
of the lemmas are all immediate from the diagram multiplication.

\proclaim{Lemma 3.2.2}
Assume $D$ has a propagating edge, $E$, connecting node $p_1$ in the north
face to node $p_2$ in the south face.

If nodes $p_1 + 1$ and $p_1 + 2$ in the north face
are connected by a (necessarily
undecorated) edge $E'$, then $U_{p_1} D$ is the $H$-admissible diagram
obtained by removing $E'$, 
disconnecting $E$ from the north face and reconnecting it to
node $p_1 + 2$ in the north face,
and installing a new undecorated edge between
points $p_1$ and $p_1 + 1$ in the north face.  The edge corresponding
to $E$ retains
its original decoration status.
\endproclaim

\proclaim{Lemma 3.2.3}
Assume that in the north face of $D$, nodes $i$ and $i+1$ are
connected by a decorated edge, $e_1$, and nodes $i+2$ and $i+3$ are connected
by an undecorated edge, $e_2$.  Assume also that $i > 1$.  Then $U_i U_{i+1}
D$ is the $H$-admissible diagram obtained from $D$ by exchanging $e_1$
and $e_2$.  This procedure has an inverse, since $D = U_{i+2} U_{i+1}
U_i U_{i+1} D$.
\endproclaim

\proclaim{Lemma 3.2.4}
Assume that in the north face of $D$, nodes 1 and 2 are connected by
a decorated edge, and nodes 3 and 4 are connected by an undecorated
edge.  Then the $H$-admissible diagram $\rho(\a) D$ is that obtained
from $D$ by decorating the edge connecting nodes 3 and 4.
\endproclaim
\proclaim{Lemma 3.2.5}
Assume that in the north face of $D$, nodes 1 and 2 are connected by
a decorated edge, $E$, and nodes 3 and 4 are connected by an undecorated
edge.  Then the $H$-admissible diagram $U_3 \rho(\z) D$ is that obtained
from $D$ by removing the decoration on $E$.
\endproclaim

\proclaim{Lemma 3.2.6}
Assume that in the north face of $D$, nodes $i$ and $i+1$ are
connected by an undecorated edge, $e_1$, and nodes $j < i$ and $k > i+1$ are
connected by an edge, $e_2$.  Assume also that $j$ and $k$ are chosen
such that $|k - j|$ is minimal.  Then $D$ is of the form $U_i D'$, where
$D'$ is an $H$-admissible diagram which 
is the same as $D$ except as regards the edges connected to nodes $j,
i, i+1, k$ in the north face.  Nodes $j$ and $i$ in $D'$ are connected
to each other by an edge with the same decoration as $e_2$, and nodes
$i+1$ and $k$ are connected to each other by an undecorated edge.
\endproclaim

As an illustration of what is going on, we present 
diagrammatic versions of these lemmas in Figure 8.  A hollow circle
indicates the site of an optional decoration.

\topinsert
\topcaption{Figure 8} Respective 
illustrations of lemmas 3.2.2 to 3.2.6 \endcaption
\centerline{
\hbox to 4.888in{
\vbox to 3.402in{\vfill
        \includegraphics{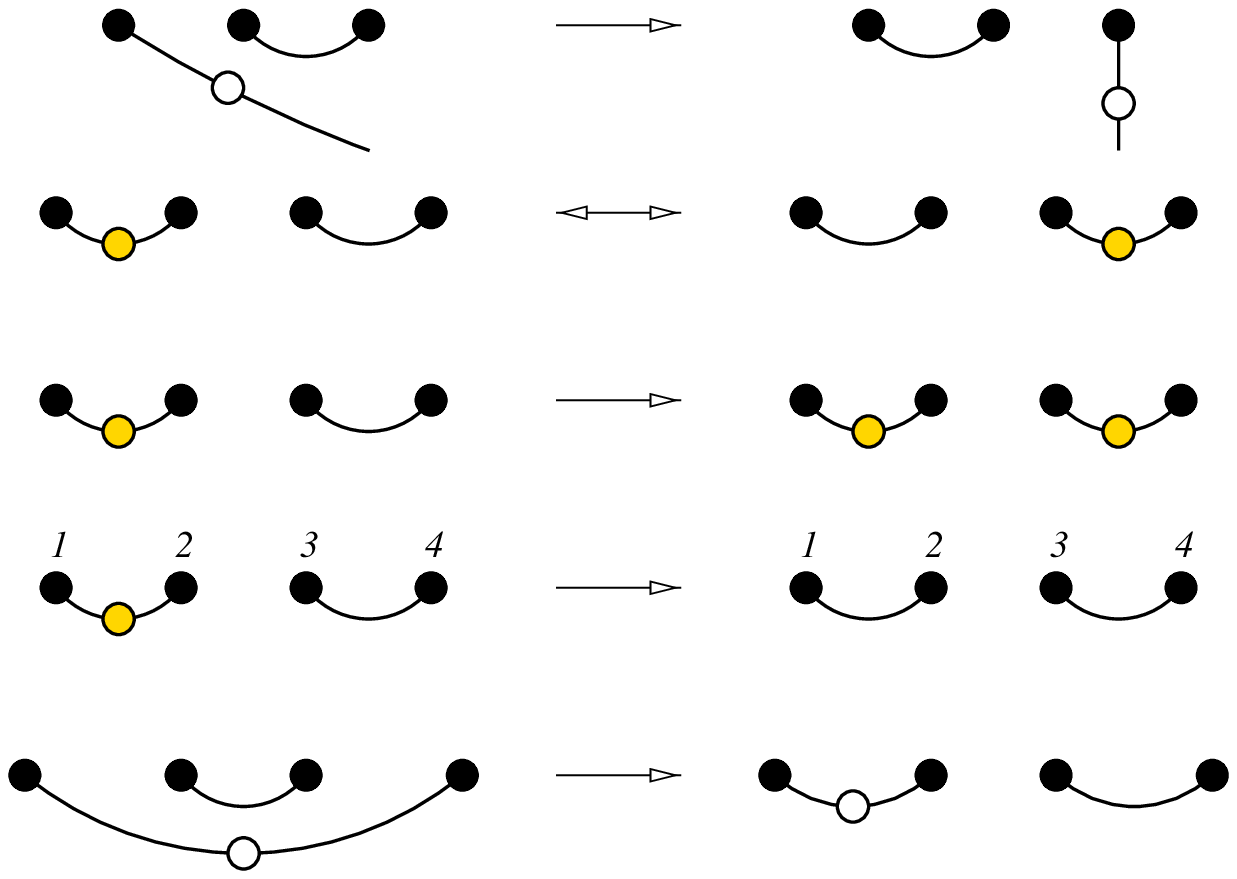}
}
\hfill}
}
\endinsert

\remark{Remark 3.2.7}
The algebra $\Delta_{n+1}$ has an anti-automorphism, $*$, which
reflects each diagram in the east-west line.  Therefore, all of the
five previous lemmas have corresponding statements about the south
faces.
\endremark

Using these five lemmas, we can prove the following result.

\proclaim{Proposition 3.2.8}
Any $H$-admissible diagram $D$ with $n+1$ strands $(n \geq 2)$
can be written as a word in the images under $\rho$ 
of the generating set $G_n$.
\endproclaim

\demo{Proof}
The case $n = 2$ is done by Lemma 3.2.1.

Iteration of Lemma 3.2.6 reduces the consideration to diagrams $D$
where all the non-propagating edges connect adjacent points.  

We restrict ourselves to the nontrivial case where $D$ has $r > 0$
non-propagating edges.

First, we assume that $D$ has a propagating edge.

We define the diagram $D_0$, depending on $D$, which is chosen
from the eight diagrams of form $$
D_0 = G U_4 U_6 \cdots U_{2r - 2} U_{2r}
,$$ where $G$ is one of the nonidentity
diagrams for the case $n = 2$ (see Lemma 3.2.1), and
$D_0$ and $D$ share the following three properties.

\item{1.}
{If $D$ has a propagating edge meeting node 1 in the north face, then so
does $D_0$.}
\item{2.}
{If $D$ has a propagating edge meeting node 1 in the south face, then so
does $D_0$.}
\item{3.}
{The leftmost propagating edges in $D$ and $D_0$ are both decorated,
or both undecorated.}

If $D$ does not have a propagating edge, we define $$
D_0 = U_1 U_3 \cdots U_n
,$$ where $n$ is necessarily odd, and $D_0$ has no propagating edges.

We claim that $D = w_1 D_0 w_2$, where $w_1$ is a word in the
generators obtained by the moves arising from lemmas
3.2.2 to 3.2.5 as stated, and $w_2$ is similar but arises from the
reflected versions of these lemmas after applying $*$ (see Remark 3.2.7).

For reasons of symmetry, we concentrate on the word $w_1$, the other
part being similar.  To do this, we show that there is a diagram $D'$
whose top half is that of $D$ and whose bottom half is that of $D_0$,
satisfying $D' = w_1 D_0$.  If $D'$ has a propagating edge, then the
leftmost one has the same decorated status as that of $D$ or $D_0$.

We start with the diagram $D_0$.  The first stage is to move the
propagating edge (if there is one) 
so that it meets the north face at the desired
point.  This is achieved by iterations of Lemma 3.2.2.

Next, we generate all the decorated, non-propagating
edges we desire, using Lemma
3.2.4.  (Note that if we have to do this, then $D_0$ has a decorated
edge connecting nodes 1 and 2 in the north face.)  
After these edges are formed, we
commute them out of the way to the east using Lemma 3.2.3.  If we
require nodes 1 and 2 in $D$ to be connected by an undecorated edge,
this can be arranged by using Lemma 3.2.5 once.  (Note that the
definition of $H$-admissible implies that $D$ must have two other points
connected by a non-decorated edge in this case, so this is possible.)
We then reach the diagram $D'$ by further applications of Lemma 3.2.3.

The proof now follows.
\qed\enddemo

\vskip 20pt

\subhead 3.3 $\Delta_{n+1}$ as a cellular algebra \endsubhead

In order to count the dimension of $\Delta_{n+1}$, it helps first to
understand its structure as a cellular algebra.  One can then compare
the sizes of the cells to those arising from $TL(H_n)$ as in [\xfanb,
Proposition 7.3.2].

It is convenient to introduce a dyadic notation for the $H$-admissible
diagrams similar to that used for $TL(A_n)$ in [\xwes, \S5].

\definition{Definition 3.3.1}
Let $D$ be an $H$-admissible diagram for $\Delta_{n+1}$.  
Remove all the propagating edges from $D$,
then take the upper half of what remains and call it $d_1$.  Invert
the lower half of $D$ in a horizontal line and call this $d_2$.  Then
$D$ may be reconstituted from the ordered pair $(d_1, d_2)$ provided
that we know whether $D$ has a decorated propagating edge or not.

We write $D = |d_1\rangle \langle d_2|$ if $D$ has no decorated
propagating edge, and $D = |d_1\rangle \langle d_2|^{\bullet}$ if $D$
has a decorated propagating edge.
\enddefinition

\proclaim{Lemma 3.3.2}
Let $R$ be an integral domain of characteristic
different from $2$, $3$ or $5$ in which the polynomial $x^2 - x - 1$
splits into distinct linear factors $(x - \g_1)(x - \g_2)$.  

Writing $\g$ for the image of $x$ in the algebra 
$\G = R[x]/\langle x^2 - x - 1 \rangle$, we have $$
(\g - \g_1)^2 = (1-2\g_1) (\g - \g_1) \ne 0
,$$ and a similar identity holds for $\g - \g_2$.
\endproclaim

\demo{Proof}
This is immediate.  Note that $1 - 2\g_1 \ne 0$ because we are not in
characteristic 5.
\qed\enddemo

\definition{Definition 3.3.3}
Let $R$ satisfy the hypotheses of Lemma 3.3.2, and let
$|d_1 \rangle \langle d_2|$ be an $H$-admissible diagram.  Then we define $$
|d_1 \rangle \langle d_2|_1 :=
|d_1 \rangle \langle d_2|^{\bullet}
- \g_1 |d_1 \rangle \langle d_2|
$$ and $$
|d_1 \rangle \langle d_2|_2 :=
|d_1 \rangle \langle d_2|^{\bullet}
- \g_2 |d_1 \rangle \langle d_2|
.$$
\enddefinition

We recall the definition of a cellular algebra from [\xgl]:

\definition{Definition 3.3.4}  Let $R$ be a commutative
ring with identity.  A {\it cellular algebra} over $R$ is an associative unital
algebra, $A$, together with a cell datum $(\Lambda, M, C, *)$ where

\item {\rm 1.}
{$\Lambda$ is a poset.  For each $\l \in \Lambda$, $M(\l)$ is a finite set
(the set of ``tableaux'' of type $\l$) such that $$
C : \coprod_{\l \in \Lambda} \left( M(\l) \times M(\l) \right) \rightarrow A
$$ is injective with image an $R$-basis of $A$.}
\item {\rm 2.}
{If $\l \in \Lambda$ and $S, T \in M(\l)$, we write $C(S, T) = C_{S, T}^{\l}
\in A$.  Then $*$ is an $R$-linear involutory anti-automorphism 
of $A$ such that
$(C_{S, T}^{\l})^* = C_{T, S}^{\l}$.}
\item {\rm 3.}
{If $\l \in \Lambda$ and $S, T \in M(\l)$ then for all $a \in A$ we have $$
a . C_{S, T}^{\l} \equiv \sum_{S' \in M(\l)} r_a (S', S) C_{S', T}^{\l}
\mod A(< \l),
$$ where  $r_a (S', S) \in R$ is independent of $T$ and $A(< \l)$ is the
$R$-submodule of $A$ generated by the set $$
\{ C_{S'', T''}^{\mu} : \mu < \l, S'' \in M(\mu), T'' \in M(\mu) \}
.$$} \enddefinition

We now define our versions of the sets in the above definition.

Let $\Lambda$ be the set of symbols $\{0\} \cup \{1, 2, \ldots, k, 1^{\bullet},
2^{\bullet}, \ldots, k^{\bullet}\}$, where $k$ is a natural number
such that $k < (n+1)/2$, together with the symbol $(n+1)/2$ if
$n$ is odd.  We put a partial order $<$ on these symbols by declaring
that $i < j$ if $|i| > |j|$, where $|i| = i$ if $i$ is a natural
number, and $|i^{\bullet}| = i$.

If $\lambda \in \Lambda$, the set $M(\l)$ has elements parametrised by
the half-diagrams $|d_1\rangle$ arising from $H$-admissible
diagrams with $|\l|$ non-propagating edges in each half of the
diagram.

The antiautomorphism $*$ corresponds to top-bottom inversion of an
$H$-admissible diagram.

The map $C$ takes elements $d_1$ and $d_2$ from $M(\l)$ and
produces the element $C(d_1, d_2)$ which is defined to be $$
|d_1 \rangle \langle d_2|_1
$$ if $\l$ is a natural number or $$
|d_1 \rangle \langle d_2|_2
$$ otherwise, unless $\l = 0$ or $\l = (n+1)/2$, in which case $C(d_1, d_2)$ is
given by $$
|d_1 \rangle \langle d_2|
.$$  Note that the identity element appears in the image of $C$.

\proclaim{Theorem 3.3.5}
Let $R$ be a ring satisfying the hypotheses of Lemma 3.3.2.
Then the algebra $\Delta_{n+1}$ over the ring $R[v, v^{-1}]$
has a cell datum $(\Lambda, M, C, *)$, where the sets are given as above.
\endproclaim

\demo{Proof}
The proof is largely straightforward.  The fact that $\g_1$ and $\g_2$
are distinct ensures that the image of $C$ is a basis for
$\Delta_{n+1}$.  

The only other nontrivial part is the verification
of axiom 3.  Consider the product of two basis elements $B_1$ and $B_2$
parametrised by the respective elements $\l$ and $\l'$ of $\Lambda$.
The only difficulty arises when $0 < \l, \l' < {{n+1}\over 2}$,
so we concentrate on this case.
It is convenient to think of each of the diagrams $B_1$ and $B_2$ as
having a propagating edge
decorated by one of the elements $\g - \g_1$ or $\g - \g_2$ of $\G$,
where an ordinary decorated edge is thought of as being decorated by
$\g \in \G$, and an undecorated one as being decorated by $1 \in \G$.
(Note that the third relation in Figure 5 corresponds to the equation
$\g^2 = \g + 1$.)

We will assume that $|\l'| \geq |\l|$; the other case is similar.
Let us define $\g_i$ by saying that the propagating edge of
$B_2$ carries the element $\g - \g_i$ (where $i \in \{1, 2\}$).
If the product $B = B_1 B_2$ is a tangle with strictly 
fewer propagating edges than $B_2$ (and therefore fewer than $B_1$,
since $|\l'| \geq |\l|$)
then it is clear that $B$ is a linear combination of
basis elements corresponding to elements $r \in \Lambda$ with $r <
\l'$, so axiom 3 holds.

The other possibility is that the
product $B$ has the same number of propagating edges as
$B_2$ and, furthermore, that the leftmost propagating edge, $E$, of
$B$ contains (as a segment) the leftmost propagating edge of $B_2$.  
The edge $E$ therefore carries the generalized decoration $\g - \g_i$, and
possibly other decorations of various kinds.

Lemma 3.3.2 shows that if we multiply together all the decorations on
the edge $E$ (where an ordinary decoration corresponds to $\g$, as
before), we obtain a (possibly zero) multiple of $\g - \g_i$.  If
we obtain zero then the product $B$ is zero and there is nothing more
to prove.  Otherwise, $B$ is a linear combination of basis elements
whose leftmost propagating edges all carry $\g - \g_i$, namely a
combination corresponding to the element $\l' \in \Lambda$.  Since the
structure constants are not affected by the pattern of non-propagating edges in
the south face of $B_2$, axiom 3 follows.
\qed\enddemo

\vskip 20pt

\subhead 3.4 Faithfulness of the diagram representation \endsubhead

Using the results of \S3.3, we can enumerate the number of
$H$-admissible diagrams of various types.

\proclaim{Lemma 3.4.1}
The size of the set $M(\l)$ associated with the algebra
$\Delta_{n+1}$ is equal to $$
{{n + 1} \choose |\l|} - 1
,$$ unless $|\l| = 0$ in which case the set has size 1.
\endproclaim

\demo{Proof}
If we generalized the $H$-admissible diagrams by excluding parts (i)
and (ii) of Definition 2.2.1, then it would follow from [\xmsa,
Proposition 2] that there would be ${{n+1} \choose k}$ half diagrams
with $k$ non-propagating edges.  If $k > 0$ then the force of part (ii) of
Definition 2.2.1 is to exclude just one element: the one with an
undecorated edge connecting points 1 and 2 and a decorated edge
connecting points $2m + 1$ and $2m + 2$ for $m < k$ (see Figure 4).  
This proves the
assertion for $|\l| > 0$, and the assertion for $|\l| = 0$ is trivial.
\qed\enddemo

\proclaim{Theorem 3.4.2}
Working over $\A$, the ranks of $\Delta_{n+1}$ and $TL(H_n)$ are identical.
Therefore, $\rho$ is an isomorphism.
\endproclaim

\demo{Proof}
Stembridge [\xsteb, \S3.4] proves that the number of ``fully
commutative'' elements in a Coxeter group of type $H_n$ is given by $$
1 + \sum_{\l \in \Lambda, |\l| > 0} \left( {{n + 1} \choose |\l|} - 1
\right) ^2 = {{2n+2} \choose {n+1}} - 2^{n+2} + n + 3
.$$  Graham [\xgra, Theorem 6.2] shows that these fully commutative
elements index a basis for $TL(H_n)$.  It follows from Lemma 3.4.1
that the rank of $\Delta_{n+1}$ is the same as the rank of $TL(H_n)$.
The fact that $\rho$ is an isomorphism follows from
Proposition 3.1.2 and Proposition 3.2.8.
\qed\enddemo

\remark{Remark 3.4.3}
The fact that $TL(H_n)$ is cellular if $x^2 - x - 1$ splits has been
observed by Graham [\xgra, Remark 9.8], although a cell datum is not
explicitly given.
\endremark

\remark{Remark 3.4.4}
The similarity with the blob algebra of [\xmsa] which is touched upon
in the proof of Lemma 3.4.1 goes further.  If conditions (i) and (ii)
of Definition 2.2.1 are dropped, then the resulting algebra is
isomorphic to the algebra of [\xmsa], although the isomorphism is not
canonical.
\endremark

\vskip 20pt

\head 4. Applications \endhead

We now examine some applications of theorems 3.3.5 and 3.4.2.  

\subhead 4.1 Positivity Properties \endsubhead

We have
seen how the diagram basis for $TL(H_n)$ can be expressed as monomials
in a certain set of algebra generators (this follows from Proposition
3.2.8 and Theorem 3.4.2).  We now show that this diagram basis has a
positivity property.

\proclaim{Proposition 4.1.1}
Assume we are working over the ring $\A$.

Any basis element occurring with nonzero coefficient
in the product of two basis elements $D_1 D_2$ associated with the
respective elements $\l_1$ and $\l_2$ of $\Lambda$ occurs with
coefficient $c [2]^k$, where $c$ is a positive integer and $k \leq
\max(|\l_1|, |\l_2|)$.

In particular, the
structure constants of the basis of diagrams are polynomials in ${\Bbb
N}[v, v^{-1}]$.
\endproclaim

\demo{Proof}
Note that the simplification process given by the rules in Figure 5
preserves positivity, and that $[2]$ is a polynomial in ${\Bbb N}[v, v^{-1}]$.

It is immediate that there cannot be more loops forming in the diagram
multiplication than there were non-propagating edges in each half of either
$D_1$ or $D_2$, which proves the assertion about the number $k$.
\qed\enddemo

Note that any basis obtained from monomials in the original set of
generators cannot have this property: consider the monomial $E_1 E_2
E_1 E_2 E_1$.

\vskip 20pt

\subhead 4.2 Semisimplicity \endsubhead

We recall from the theory of cellular algebras in [\xgl, \S2] that there is
a bilinear form $\phi_{\l}(d_1, d_2) = \langle d_1, d_2 \rangle$
on the cell module $W(\l)$.  (Recall that the module $W(\l)$ has a
basis parametrised by the elements of $M(\l)$ for a fixed $\l$.)
This form is defined from the equation $$
C(e_1, d_1) C(d_2, e_2)
= \langle d_1, d_2 \rangle C(e_1, e_2) \mod A(< \l)
.$$  This is independent of the choice of $e_1$ and $e_2$, where $e_1,
e_2, d_1, d_2$ are all elements of $M(\l)$ for the same $\l$.

The following result proves [\xfanb, Conjecture 7.3.1].

\proclaim{Theorem 4.2.1}
Let $R$ be a field satisfying the hypotheses of Lemma 3.3.2.  Then
the algebra $TL(H_n)$ over $R$ is semisimple, and all the cell modules $W(\l)$
are irreducible and pairwise inequivalent.
\endproclaim

\demo{Proof}
It is enough by [\xgl, Theorem 3.8] to prove that $\phi_{\l}$ is
nondegenerate for each $\l$.

Choose an element $\l \in \Lambda$ and two elements $d_1, d_2 \in
M(\l)$, where possibly $d_1 = d_2$.  Now consider 
$v^{-|\l|} \phi_{\l}(d_1, d_2) = v^{-|\l|} \langle d_1, d_2
\rangle$.

It follows from Proposition 4.1.1 that $v^{- |\l|} \langle d_1, d_2
\rangle$ is a polynomial in $v^{-1}$; furthermore, the constant term
of the polynomial is zero unless
$d_1 = d_2$.  To see this, we use the fact that loops
carrying a single blob result in annihilation of the associated
diagram, and loops which carry 0 or 2 blobs both correspond to
multiplication by $[2]$.
In the case where $d_1 = d_2$, Lemma 3.3.2 shows that the constant
term of the polynomial is nonzero and equal to $(1 - 2\g_1)$ or
$(1 - 2\g_2)$, depending on the $\l$ which we are considering.

We have now constructed an almost orthogonal basis
(\idest orthogonal modulo the span of strictly negative powers of $v$)
for the module $W(\l)$ with respect to $\phi_{\l}$.  
It follows that $\phi_{\l}$ is nondegenerate, as required.
\qed\enddemo

Note that Lemma 3.4.1 now tells us the dimensions of the irreducible
modules.  This confirms [\xfanb, Conjecture 7.3.3].

\vskip 20pt

\subhead 4.3 Branching rules \endsubhead

In this section, we continue to assume that the base ring of $TL(H_n)$
is a field in which $x^2 - x - 1$ splits into distinct linear factors,
although we do not assume that $TL(H_n)$ is semisimple.
The diagram calculus we have developed allows us to study the
behaviour of the cell modules $W(\l)$ for $TL(H_n)$ upon restriction
to $TL(H_{n-1})$ (assuming $n$ is at least 3).  
The embedding of $TL(H_{n-1})$ into $TL(H_n)$ 
is the natural one arising from the identification of the algebra
generators, or the addition of a vertical edge on the east of the
diagram.  To describe the branching rules, it is convenient to make the
following definition.

\definition{Definition 4.3.1}
Let $\l \in \Lambda = \Lambda(H_n)$, and suppose that 
$0 < |\l| < {{n+1}\over 2}$.  
We define an element $\l - 1 \in \Lambda(H_{n-1})$  
as follows: $$\l - 1 := \cases
0 & \text{\rm if $|\lambda| = 1$,}\cr
(i - 1)^{\bullet} & \text{\rm if $|\lambda| \ne 1$,
$\lambda = i^{\bullet}$ where $i \in {\Bbb N}$,
and $i - 1 \ne n/2,$}\cr
i - 1 & \text{\rm otherwise.}\cr
\endcases$$
\enddefinition

\proclaim{Proposition 4.3.2}
Let $W(\l, n)$ be a cell module for $TL(H_n)$.  Then, after restriction
to $TL(H_{n-1})$, $W(\l, n)_{n-1}$ 
has a filtration by the cell modules
$W(\l', n-1)$ of $TL(H_{n-1})$ described as follows.

If $\l = 0$ then restriction gives the trivial module corresponding to the
poset element $0$.

If $|\l| = 1$ then the composition factors occurring correpond to the
poset elements $0$ and $\l$, each with multiplicity $1$.

If $\l = {{n+1}\over 2}$ then the composition factors occurring
correspond to the poset elements $0$, ${{n-1}\over 2}$ and $\left({{n-1}\over
2}\right)^{\bullet}$, each with multiplicity $1$.

For other values of $\l$, the composition factors occurring correspond
to the poset elements $0$, $\l - 1$ and $\l$, each with multiplicity $1$.
\endproclaim

\demo{Proof}
We first tackle the fourth case, dealing with the general value of
$\l$.  

The key observation, which is familiar from the diagram calculi
of other types, is as follows.  The half-diagrams in $M(\l)$ which
have a non-connected point at the eastern extreme form a submodule for
$TL(H_{n-1})$ on restriction, corresponding to removal of the
easternmost point.  This is canonically isomorphic to the module
corresponding to $\l$ in $\Lambda(H_{n-1})$.

The quotient module associated with this submodule is obtained by taking the
other elements of $M(\l)$ and, for each one,
removing the easternmost point and the edge connected to it.  However,
this is not the same as one of the cell modules for $TL(H_{n-1})$, because
an inadmissible half-diagram occurs.  (All the edges in this are
decorated, except for the one connecting points $1$ and $2$.)
The reason that this arises is that in the
original diagram, the two easternmost points could have been connected by an
undecorated edge, making the half-diagram admissible, but this edge
was removed in the restriction process.  The admissible diagrams
arising from this procedure span a cell module isomorphic to that
parametrised by $\l - 1$.  The appearance of the inadmissible diagram
corresponds to a top quotient isomorphic to the trivial module.

The cases $|\l| = 0$ and $|\l| = 1$ can be obtained by degenerate
versions of this technique.  

The case $\l = {{n+1}\over 2}$ is the most subtle.  To deal with it,
it is convenient to modify the basis for the cell module $W(\l, n)$.
One of the half-diagrams in $M(\l)$ becomes inadmissible once the
rightmost point and its associated edge have been removed; we denote
this half-diagram by $d_0$.  

The other half-diagrams fall naturally into pairs as we now describe.
First, note that any edge of a half-diagram of $M(\l)$ (for $\l
= {{n+1}\over 2}$) is exposed to the west face, since there are no
propagating edges involved.  There is therefore an involution on
$M(\l) \backslash \{d_0\}$ given by changing the decorated status of
the edge connected to the rightmost point.  The orbits are all of size 2.  

If $d$ and $d^{\bullet}$ are two elements in the
same orbit (where $d^{\bullet}$ carries the extra decoration), we define
new basis elements $d_1$ and $d_2$ for $W(\l, n)$ by $$
d_1 := d^{\bullet} - \g_1 d
$$ and $$
d_2 := d^{\bullet} - \g_2 d
.$$  

Analysis of the diagrams now shows that the diagram $d_0$ corresponds
to a top quotient of $W(\l, n)$
isomorphic to the trivial module.  Furthermore, the
submodule of $W(\l, n)$ spanned by the new basis elements $d_i$ breaks up as a
direct sum: the span of the elements $d_1$ is canonically 
isomorphic to $W\left({{n-1}\over 2}, n-1\right),$ and
the span of the elements $d_2$ is canonically isomorphic to $
W\left(\left({{n-1}\over 2}\right)^{\bullet}, n-1\right).$
\qed\enddemo

\vskip 20pt

\head Acknowledgements \endhead

The author is grateful to C.K. Fan for some useful discussions, and
to the referee for some helpful comments.
\vfill\eject

\vskip 1cm

\head References \endhead

\item{[\xfanb]}
{C.K. Fan, {\it Structure of a Hecke Algebra Quotient}, Jour.
Amer. Math. Soc. {\bf 10} (1997), 139--167.}
\item{[\xfga]}
{C.K. Fan and R.M. Green, {\it On the affine Temperley--Lieb
algebras}, Jour. L.M.S., to appear.}
\item{[\xfy]}
{P.J. Freyd and D.N. Yetter,
{\it Braided compact closed categories with applications to low
dimensional topology}, Adv. Math. {\bf 77} (1989), 156--182}
\item{[\xgra]}
{J.J. Graham, Ph.D. Thesis, University of Sydney, 1995.}
\item{[\xgl]}
{J.J. Graham and G.I. Lehrer, {\it Cellular Algebras}, Invent. Math. {\bf 123}
(1996), 1--34.}
\item{[\xrmgl]}
{R.M. Green, {\it Generalized Temperley--Lieb algebras and decorated
tangles}, Jour. Knot Th. Ram., to appear.}
\item{[\xhum]}
{J.E. Humphreys,
{\it Reflection Groups and Coxeter Groups}, Cambridge University
Press, 1990.}
\item{[\xjon]}
{V.F.R. Jones, {\it Hecke algebra representations of braid groups and
link polynomials}, Ann. of Math. (2) {\bf 126} (1987), 335--388.}
\item{[\xk]}
{L.H. Kauffman,
{\it State models and the Jones polynomial},
Topology {\bf 26} (1987), 395--407}
\item{[\xmsa]}
{P. Martin and H. Saleur, {\it The blob algebra and the periodic
Temperley--Lieb algebra}, Lett. Math. Phys., {\bf 30} (1994) no. 3, 189--206.}\
\item{[\xsteb]}
{J.R. Stembridge,
{\it The Enumeration of Fully Commutative Elements of Coxeter Groups},
J. Alg. Comb., to appear.}
\item{[\xtl]}
{H.N.V. Temperley and E.H. Lieb, {\it Relations between percolation
and colouring problems and other graph theoretical problems associated
with regular planar lattices: some exact results for the percolation problem},
Proc. Roy. Soc. London Ser. A {\bf 322} (1971), 251--280.}
\item{[\xwes]}
{B.W. Westbury, {\it The representation theory of the Temperley--Lieb
Algebras}, \newline Math. Z. {\bf 219} (1995), 539--565.}

\vfill\eject
\end